\documentclass[12pt]{amsart}
\usepackage{amssymb}
\usepackage{amsthm}
\usepackage{enumerate}
\usepackage[dvips]{graphics,color}

\theoremstyle{plain}
\newtheorem{Thm}{Theorem}[section]
\newtheorem{Lem}[Thm]{Lemma}

\newtheorem{Cor}[Thm]{Corollary}

\theoremstyle{remark}

\small\normalsize
\numberwithin{equation}{section}
\def\beginpf{\noindent {\bf Proof:} \quad}
\newcommand{\boite}{\mbox{} \hfill $\square$}
\def\endpf{\boite\par\medskip}

\def\CC{{\mathbb C}}

\def\NN{{\mathbb N}}
\def\RR{{\mathbb R}}
\def\DD{{\mathbb D}}
\def\TT{{\mathbb T}}

\def\f1{\mathbb{1}}
\def\HH{{\mathcal H}}

\def\la{\lambda}
\def\H{\mathcal H}

\def\dfrac{\displaystyle\frac}

\begin{document}

\title{Boundary behavior of functions in the de Branges--Rovnyak spaces}

\author{Emmanuel Fricain, Javad Mashreghi}

\address{Universit\'e de Lyon; Universit\'e Lyon 1; Institut Camille Jordan CNRS UMR 5208; 43, boulevard du 11 Novembre 1918, F-69622 Villeurbanne}
\email{fricain@math.univ-lyon1.fr}

\address{D\'epartement de math\'ematiques et de statistique,
         Universit\'e Laval,
         Qu\'ebec, QC,
         Canada G1K 7P4.}
\email{Javad.Mashreghi@mat.ulaval.ca}

\thanks{This work was supported by NSERC (Canada) and FQRNT (Qu\'ebec). A part of this work was done while the first author was visiting McGill University. He would like to thank this institution for its warm hospitality.}



\begin{abstract}
This paper deals with the boundary behavior of functions in the de
Branges--Rovnyak spaces. First, we give a criterion for the
existence of radial limits for the derivatives of functions in the
de Branges--Rovnyak spaces. This criterion generalizes a result of
Ahern-Clark. Then we prove that the continuity of all functions in a
de Branges--Rovnyak space on an open arc $I$  of the boundary is
enough to ensure the analyticity of these functions on $I$. We use
this property in a question related to Bernstein's inequality.
\end{abstract}

\maketitle

\section{Introduction}
For $0<p\leq\infty$, let $H^p(\DD)$ denote the classical Hardy space
of analytic functions on the unit disc $\DD:=\{z\in\CC:|z|<1\}$. As
usual, we also treat $H^p(\DD)$ as a closed subspace of
$L^p(\TT,m)$, where $\TT:=\partial\DD$ and $m$ is the normalized arc
length measure on $\TT$. Let $b$ be in the unit ball of
$H^\infty(\DD)$. Then the canonical factorization of $b$ is $b=BF$,
where
\[
B(z)= \gamma \prod_{n}
\dfrac{|a_n|}{a_n}\frac{a_n-z}{1-\overline{a_n}z},\qquad (z\in\DD),
\]
is the Blaschke product with zeros $a_n\in\DD$ satisfying the
Blaschke condition $\sum_n(1-|a_n|)<~+\infty$, $\gamma$ is a
constant of modulus one, and $F$ is of the form
\[
F(z)=\exp\left(-\int_\TT\frac{\zeta+z}{\zeta-z}\,d\sigma(\zeta)
\right),\qquad (z\in\DD),
\]
where $d\sigma=-\log|b|\,dm+d\mu$ and $d\mu$ is a positive singular
measure on $\TT$. In the definition of $B$, we assume that
$|a_n|/a_n=1$ whenever $a_n=0$. In this paper, we study some aspects
of the de Branges--Rovnyak spaces
\[
\HH(b):=(Id-T_bT_{\overline b})^{1/2}H^2.
\]
Here $T_\varphi$ denotes the Toeplitz operator defined on $H^2$ by
$T_\varphi(f)=P_+(\varphi f)$, where $P_+$ is the (Riesz) orthogonal
projection of $L^2(\TT)$ onto $H^2$. In general, $\HH(b)$ is not
closed with respect to the norm of $H^2(\DD)$. However, it is a
Hilbert space when equipped with the inner product
$$\langle\, (Id-T_b T_{\overline b})^{1/2}f, \,(Id-T_b T_{\overline b})^{1/2}g \,\rangle_b=\langle f,g
\rangle_2,$$ where $f$ and $g$ are chosen so that $$f,g \, \bot \,
\hbox{ker }(Id-T_b T_{\overline b})^{1/2}.$$ As a very special case,
if $|b|=1$ a.e. on $\TT$, or equivalently when $b$ is an inner
function for the unit disc, then  $Id-T_bT_{\overline b}$ is an
orthogonal projection and the $\HH(b)$ norm coincides with the $H^2$
norm. In this case, $\HH(b)$ becomes a closed (ordinary) subspace of
$H^2(\DD)$, which coincides with the shift-coinvariant subspace
$K_b:=H^2\ominus bH^2$.

This paper treats two questions related to the boundary behavior of
functions in $\HH(b)$. The first of these concerns the existence of
radial limits for the derivatives of functions in the de
Branges--Rovnyak spaces. More precisely, given a non-negative
integer $N$, we are interested in finding a characterization of
points $\zeta_0\in\TT$ such that every function $f$ in $\HH(b)$ and
its derivatives up to order $N$ have radial limits at $\zeta_0$.
Ahern and Clark \cite{ak70} studied this question when $b$ is an
inner function and they got a characterization in terms of the zeros
sets $(a_n)$ and the measure $\mu$. In Section 3, we show that their
methods in \cite{ak70,ak71} can be extended in order to obtain
similar results for the general de Branges--Rovnyak spaces $\HH(b)$,
where $b$ is an arbitrary element of the unit ball of $H^\infty$.
Let us also mention that Sarason \cite[page 58]{sarason} has
obtained another criterion in terms of the measure whose Poisson
integral is the real part of $\dfrac {\la+b}{\la-b}$, with
$\la\in\TT$. Recently, Bolotnikov and Kheifets \cite{Bolotnikov}
gave a result, in some sense more algebraic, in terms of the
Schwarz-Pick  matrix.

Our second theme is related to the analytic continuation of
functions in $\HH(b)$ through a given open arc of $\TT$. In
\cite{helson}, in the case where $b$ is an inner function, Helson
proved that every function in $K_b$ has an analytic continuation
through an open arc $I$ of $\TT$ if and only if $b$ has an analytic
continuation through $I$. Then, in \cite[page 42]{sarason}, Sarason
extended this result to the de Branges--Rovnyak spaces $\HH(b)$,
when $b$ is an extreme point of the unit ball of $H^\infty$. In the
last section, we study the question of continuity on the open arc
$I$ for functions in $\HH(b)$. In particular, we show that the
continuity on some open arc of the boundary of all functions in
$\mathcal H(b)$ implies the analyticity on this arc. We apply this
remarkable property to discuss a possible generalization of the
Bernstein's inequality obtained by Dyakonov \cite{Dyak02} in the
model space $K_b$.

\section{Preliminaries}
We first recall some basic well-known facts concerning reproducing
kernels in $\HH(b)$. For any $\la\in\DD$,  the linear functional
$f\longmapsto f(\la)$ is bounded on $H^2(\DD)$ and thus, by Riesz'
theorem, it is induced by a unique element $k_\la$ of $H^2(\DD)$. On
the other hand, by Cauchy's formula, we have
$$f(\la)=\frac 1{2\pi}\int_0^{2\pi}\frac {f(e^{i\vartheta})}{1-\la e^{-i\vartheta}}\,d\vartheta,\qquad (f\in H^2(\DD),\,\la\in\DD),$$
and thus
$$k_\la(z)=\frac 1{1-\overline \la z},\qquad (z\in\DD).$$
Now, since $\HH(b)$ is contained contractively in $H^2(\DD)$, the
restriction to $\HH(\DD)$ of the  evaluation functional at
$\la\in\DD$ is  a bounded linear functional on $\HH(\DD)$. Hence,
relative to the inner product in $\HH(b)$, it is  induced by a
vector $k_\la^b$ in $\HH(b)$. In other words, for all $f\in\HH(b)$,
we have
$$f(\la) = \langle f,k_\la^b \rangle_b.$$
But if $f=(Id-T_b T_{\overline b})^{1/2}f_1\in\HH(b)$, we have
$$\langle f,(Id-T_b T_{\overline b})k_\la\rangle_b=\langle f_1,(Id-T_b T_{\overline b})^{1/2}k_\la\rangle_2=\langle f,k_\la \rangle_2=f(\la),$$
which implies that
$$k_\la^b=(Id-T_b T_{\overline b})k_\la.$$
Finally, using the well known result $T_{\overline
b}k_w=\overline{b(w)}k_w$, we obtain
$$k_\la^b(z)=\frac{1-\overline{b(\la)}b(z)}{1-\overline\la z},\qquad (z\in\DD).$$

We know (see \cite[page 11]{sarason})  that $\HH(b)$ is invariant
under the backward shift operator $S^*$ and, in the following, we
use extensively the contraction $X:=S^*|\HH(b)$. Its adjoint satisfies
the important formula
\begin{eqnarray}\label{eq:operateurX}
X^*h=Sh-\langle h,S^*b\rangle_b \, b,
\end{eqnarray}
for all $h\in\HH(b)$ (see \cite[page 12]{sarason}).

We end this section by recalling the definition of the spectrum of a
function $b$ in the unit ball of $H^\infty(\DD)$ (see \cite[page
103]{nikolski-controle2}). A point $\lambda\in\overline\DD$ is said
to be regular (for $b$) if either $\lambda\in\DD$ and
$b(\lambda)\not=0$, or $\lambda\in\TT$ and $b$ admits an analytic
continuation across a neighbourhood
$V_\lambda=\{z:|z-\lambda|<\varepsilon\}$ of $\lambda$ with $|b|=1$
on $V_\lambda \cap \TT$. The spectrum of $b$, denoted by
$\sigma(b)$, is then defined as the complement in $\overline\DD$ of
all regular points of $b$.

\section{Existence of derivatives for functions of de Branges--Rovnyak spaces}

We first begin with a lemma which is essentially due to Ahern-Clark
\cite[Lemma 2.1]{ak70}.

\begin{Lem}\label{lem:ac}
Let $S_1,\dots,S_p$ be bounded commuting operators of norm less or
equal to 1 on a Hilbert space $X$. Let
$(\lambda_1,\dots,\lambda_p)\in\TT^p$ such that $Id-\lambda_jS_j$ is
one to one. Furthermore, let
$(\lambda_1^{(n)},\dots,\lambda_p^{(n)})\in\DD^p$ tends
nontangentially to $(\lambda_1,\dots,\lambda_p)$ as $n\to +\infty$.
Then, for any $y\in X$, the sequence
$w_n:=(Id-\lambda_1^{(n)}S_1)^{-1}\dots(Id-\lambda_p^{(n)}S_p)^{-1}y$
is uniformly bounded if and only if $y$ belongs to the range of the
operator $(Id-\lambda_1S_1)\dots (Id-\lambda_p S_p)$, in which case,
$w_n$ tends weakly to $w_0:=(Id-\lambda_1
S_1)^{-1}\dots(Id-\lambda_p S_p)^{-1}y$.
\end{Lem}

\beginpf
If $\|S_j\|<1$, then the operator $Id-\lambda_jS_j$ is invertible and $(Id-\lambda_j^{(n)}S_j)^{-1}$ tends to $(Id-\lambda_jS_j)^{-1}$ in operator norm, as $n\to +\infty$. Therefore, we see that we can assume that all operators $S_j$ are of norm equal to $1$. This case is precisely the result of Ahern-Clark.

\endpf

The following result gives a criterion for the existence of the
derivatives for functions of $\HH(b)$ and it generalizes the
Ahern-Clark result.

\begin{Thm}\label{Thm:derivative-DD}
Let $b$ be a point in the unit ball of $H^\infty(\DD)$ and let
{\small\begin{equation}\label{eq:factorisation-disque}
b(z)=\gamma\prod_n \left(\, \frac{|a_n|}{a_n}\frac
{a_n-z}{1-\overline{a}_nz} \,\right) \,\,%
\exp\left( -\int_\TT \frac {\zeta+z}{\zeta-z}\,d\mu(\zeta)\right)
\,\,%
\exp\left(\int_\TT \frac
{\zeta+z}{\zeta-z}\log|b(\zeta)|\,dm(\zeta)\right)
\end{equation}}
be its canonical factorization. Let $\zeta_0\in\TT$ and let
 $N$ be a non-negative integer. Then the following are equivalent.
\begin{enumerate}
\item[$\mathrm{(i)}$] for every function $f\in\HH(b)$, $f(z),f'(z),\dots,f^{(N)}(z)$ have finite limits as $z$ tends radially to $\zeta_0$;
\item[$\mathrm{(ii)}$] for every function $f\in\HH(b)$,  $|f^{(N)}(z)|$ remains bounded as $z$ tends radially~to~$\zeta_0$;
\item[$\mathrm{(iii)}$]  $\left\|\partial^N k_z^{b}/\partial {\overline z^N} \right\|_{b}$ is bounded as $z$ tends radially to $\zeta_0$;
\item[$\mathrm{(iv)}$] ${X^*}^Nk_0^{b}$ belongs to the range of $(Id-\overline{\zeta_0}X^*)^{N+1}$;
\item[$\mathrm{(v)}$] we have
\[
\sum_n\frac {1-|a_n|^2}{|\zeta_0-a_n|^{2N+2}}+\int_0^{2\pi} \frac {d\mu(e^{it})}{|\zeta_0-e^{it}|^{2N+2}}+\int_0^{2\pi} \frac {\big| \log|b(e^{it})| \big|}{|\zeta_0-e^{it}|^{2N+2}} \,\, dm(e^{it})<+\infty.
\]
\end{enumerate}
\end{Thm}

\beginpf

$\mathrm{(i)}\Longrightarrow\mathrm{(ii)}$: it is obvious.

$\mathrm{(ii)}\Longrightarrow\mathrm{(iii)}$:  for a point $z$ in
$\DD$, the function $\displaystyle\frac{\partial^N k_z^{b}}{\partial
{\overline z}^N}$ is easily seen to be the kernel function in
$\HH(b)$ for the functional of evaluation of the $N$th derivative at
$z$:
\begin{eqnarray}\label{eq:derivee}
f^{(N)}(z)=\langle f,\frac{\partial^N k_z^{b}}{\partial {\overline
z}^N}\rangle_{b},\qquad\forall f\in\HH(b).
\end{eqnarray}
Therefore, the implication
$\mathrm{(ii)}\Longrightarrow\mathrm{(iii)}$ follows from the
principle of uniform boundedness.

The equivalence of $\mathrm{(i)}$ and $\mathrm{(iii)}$ is not new
and can be found in \cite[page 58]{sarason}.

$\mathrm{(iii)}\Longrightarrow\mathrm{(iv)}$:  using the fact that
$k_z^{b}=(Id-\overline z X^*)^{-1}k_0^{b}$ (see \cite[page
42]{sarason}), we easily get
\begin{eqnarray}\label{eq:cle1}
\frac{\partial^N k_z^{b}}{\partial {\overline z}^N}=N! (Id-\overline
z X^*)^{-(N+1)}{X^*}^N k_0^{b}.
\end{eqnarray}
We know from \cite[Lemma 2.2]{ef2005} that $\sigma_p(X^*)\subset\DD$
and thus the operator $Id-\overline{\zeta_0}X^*$ is one-to-one. By
assumption, $(Id-\overline{z_n}X^*)^{-(N+1)}{X^*}^Nk_0^b$ is
uniformly bounded for any sequence $z_n\in\DD$ tending radially to
$\zeta_0$. Hence, by Lemma \ref{lem:ac}, ${X^*}^Nk_0^b$ belongs
to the range of $(Id-\overline{\zeta_0}X^*)^{N+1}$.

$\mathrm{(iv)}\Longrightarrow\mathrm{(i)}$:  using once more Lemma \ref{lem:ac} with $p=N+1$, $S_1=\dots=S_p=X^*$, $\lambda_1=\dots=\lambda_p=\overline\zeta_0$ and $y={X^*}^N k_0^b$, we see that $\mathrm{(iv)}$ implies that
$(Id-\overline{z_n}X^*)^{-(N+1)}{X^*}^Nk_0^b$ tends weakly to
$(Id-\overline{\zeta_0}X^*)^{-(N+1)}{X^*}^Nk_0^b$, for any sequence
$z_n\in\DD$ tending radially to $\zeta$. Hence (\ref{eq:derivee})
and (\ref{eq:cle1}) imply that, for every function $f$ in $\HH(b)$,
$f^{(N)}(z)$ has a finite limit as $z$ tends radially to $\zeta_0$.
Now of course, for every $0\leq j\leq N$,  $\mathrm{(iv)}$ ensures
that ${X^*}^jk_0^b$ belongs to the range of
$(Id-\overline{\zeta_0}X^*)^{j+1}$ and similar arguments show that,
for every function $f$ in $\HH(b)$, $f^{(j)}(z)$ has a finite limit
as $z$ tends radially to $\zeta_0$.

$\mathrm{(v)}\Longrightarrow\mathrm{(iii)}$:  without loss of
generality we assume that $\zeta_0=1$. Using Leibnitz' rule, by
straightforward computations
 we obtain
\begin{equation}\label{eq:derive-explicite}
k_{\omega,N}^b(z):=\frac{\partial^N k_\omega^b}{\partial\overline\omega^N}(z)=\frac{h_{\omega,N}^b(z)}{(1-\overline\omega z)^{N+1}},
\end{equation}
with
\begin{equation}\label{eq:derive2-explicite}
h_{\omega,N}^b(z)=N!z^N-b(z)
\displaystyle\sum_{j=0}^N\binom{N}{j}\overline{b^{(j)}(\omega)}(N-j)!z^{N-j}(1-\overline\omega
z)^j.
\end{equation}
Hence, by (\ref{eq:derivee}), we have
$$
\left\|\frac{\partial^Nk_\omega^b}{\partial\overline\omega^N}\right\|_b^2=(k_{\omega,N}^b)^{(N)}(\omega),
$$
and thus, we need to prove that $(k_{r,N}^b)^{(N)}(r)$ is bounded as
$r\to 1^{-}$.

But the condition $\mathrm{(v)}$ clearly implies that
\[
\sum_n\frac {1-|a_n|^2}{|\zeta_0-a_n|^{j}}+\int_0^{2\pi} \frac {d\mu(e^{it})}{|\zeta_0-e^{it}|^{j}}+\int_0^{2\pi} \frac {\big| \log|b(e^{it})| \big|}{|\zeta_0-e^{it}|^{j}} \,\, dm(e^{it})<+\infty,
\]
for $0\leq j\leq 2N+2$ and then it follows from \cite[Lemma 4]{ak71} that
$$\lim_{r\to 1^-}b^{(j)}(r) \quad \mbox{ and } \quad \lim_{R\to 1^+}b^{(j)}(R)$$exist and are equal. Here we extend the function $b$ outside the
unit disk by the formula (\ref{eq:factorisation-disque}), which
represents an analytic function  for $|z|>1$,
$z\not=1/\overline{a_n}$. We denote this function also by $b$ and it
is easily verified that it satisfies
\begin{eqnarray}\label{eq:fonctionnelle}
b(z)=\frac 1{\,\,\overline{b(1/\overline z)}\,\,},\qquad \forall
z\in\CC.
\end{eqnarray}
Therefore,  there exists $R_0>1$ such that $b$ has $2N+1$ continuous
derivatives on $[0,R_0]$. Now take $R_0^{-1}<r<1$. Noting that $b$
can have only a finite number of real zeros, we can assume that the
interval $(R_0^{-1},1)$ is free of zeros. Then straightforward computations using (\ref{eq:derive2-explicite}) and (\ref{eq:fonctionnelle}) show that $h_{r,N}^b$ and its first $N$
derivatives must vanish at $z=1/r$. Therefore we can write, for
$s\in (0,1)$,
$$\begin{aligned}
h_{r,N}^b(s)=&\int_0^1\frac {d}{dt}h_{r,N}^b\left(\frac 1r+t(s-\frac 1r)\right)\,dt \\
=& \left(s-\frac 1r\right)\int_0^1 (h_{r,N}^b)'\left(\frac 1r +t(s-\frac 1r)\right)\,dt\\
=&\left(s-\frac 1r\right)^2\int_0^1\!\!\int_0^1 (h_{r,N}^b)''\left(\frac 1r +tu(s-\frac 1r)\right)t\,du\,dt.\\
\end{aligned}$$
Continuing this procedure, we get
\[
h_{r,N}^b(s)=\left(s-\frac
1r\right)^{N+1}\!\int_0^1\!\!\int_0^1\dots\int_0^1
(h_{r,N}^b)^{(N+1)}\left(\frac 1r +t_1t_2\dots t_{N+1}(s-\frac
1r)\right)m(t)\,dt_1\,\dots\,dt_{N+1},
\]
where $m(t)$ is a monomial in $t_1,\dots,t_{N+1}$. Hence, using (\ref{eq:derive-explicite}), we obtain
\[
k_{r,N}^b(s)=\dfrac{1}{r^{N+1}}\int_0^1\!\!\int_0^1\dots\int_0^1 (h_{r,N}^b)^{(N+1)}\left(\frac 1r +t_1t_2\dots t_{N+1}(s-\frac 1r)\right)m(t)\,dt_1\,\dots\,dt_{N+1}.
\]
But, thanks to properties of $b$, we can differentiate under the
integral sign to get
$$(k_{r,N}^b)^{(N)}(s)=\dfrac{1}{r^{N+1}}\int_0^1\!\!\int_0^1\dots\int_0^1 (h_{r,N}^b)^{(2N+1)}\left(\frac 1r +t_1t_2\dots t_{N+1}(s-\frac 1r)\right)v(t)\,dt_1\,\dots\,dt_{N+1},$$
where $v(t)$ is a monomial in $t_1,\dots,t_{N+1}$. Since
$(h_{r,N}^b)^{(2N+1)}$ is bounded on $(0,R_0)$, we deduce that
$|(k_{r,N}^b)^{(N)}(r)|\leq \frac 1{r^{N+1}}
\|(h_{r,N}^b)^{(2N+1)}\|_\infty$, which is bounded as $r\to 1^-$.

$\mathrm{(iii)}\Longrightarrow\mathrm{(v)}$: here we also assume
that $\zeta_0=1$. According to \cite[Lemma 4.2]{ak70} we can take a
sequence $(B_j)_{j\geq 1}$ of Blaschke products converging uniformly
to $b$ on compact subsets of $\DD$ and such that
$$\sum_{k} \frac {1-|a_{j,k}|^2}{|1-r a_{j,k}|^{2N+2}}\,\,\,\,{}_{\substack {\longrightarrow \\
                      j\to +\infty}}\sum_{k} \frac {1-|a_k|^2}{|1-r
a_k|^{2N+2}}+\int_0^{2\pi}\frac {d\mu(e^{it})}{|e^{it}-r|^{2N+2}}+\int_0^{2\pi}\frac {|\log|b(e^{it})||}{|e^{it}-r|^{2N+2}}\,dm(e^{it}),$$
where $(a_{j,k})_{k\geq 1}$ is the sequence of zeros of $B_j$. As
before, let $k_{\omega,N}^b:=\dfrac {\partial^N
k_\omega^{b}}{\partial\overline\omega^N}$ and let
$k_{\omega,N}^{B_j}:=\dfrac {\partial^N
k_\omega^{B_j}}{\partial\overline\omega^N}$. Hence, we have
\begin{eqnarray}\label{eq:noyau-derive}
k_{\omega,N}^{B_j}(z)=\frac{N!z^N-B_j(z) \displaystyle\sum_{p=0}^N\binom{N}{p}\overline{B_j^{(p)}(\omega)}(N-p)!z^{N-p}(1-\overline\omega z)^p}{(1-\overline\omega z)^{N+1}}
\end{eqnarray}
and thus $k_{\omega,N}^{B_j}$ tends to $k_{\omega,N}^b$ uniformly on
compact subsets of $\DD$. Therefore,
$$\lim_{j\to +\infty}(k_{\omega,N}^{B_j})^{(N)}(\omega)=(k_{\omega,N}^b)^{(N)}(\omega).$$
But,
$$\left\|\frac {\partial^N k_\omega^{b}}{\partial\overline\omega^N}\right\|_{b}^2=(k_{\omega,N}^b)^{(N)}(\omega),\qquad\hbox{and}\qquad\left\|\frac {\partial^N k_\omega^{B_j}}{\partial\overline\omega^N}\right\|_2^2=(k_{\omega,N}^{B_j})^{(N)}(\omega),$$
and condition $\mathrm{(iii)}$ implies that there exists $C_1>0$ such that, for all $0<r<1$, we have $|(k_{r,N}^b)^{(N)}(r)|\leq C_1$. Therefore, for all $0<r<1$, there exists $j_r\in\NN$, such that
for $j\geq j_r$, we have
$$\left\|\frac{\partial^N k_r^{B_j}}{\partial r^N}\right\|_2^2=|(k_{r,N}^{B_j})^{(N)}(r)|\leq C_1+1.$$
Moreover, using (\ref{eq:noyau-derive}), we see that
$$(1-rz)^{N+1} \frac{\partial^N k_r^{B_j}}{\partial r^N}(z) =N!z^N-B_j(z)g_j(z),$$
where $g_j\in H^2$. Hence, it follows from \cite[Theorem 3.1]{ak70}
that there is a constant $K$ (independent of $r$) such that
\[
\sum_{k} \frac {1-|a_{j,k}|^2}{|1-r a_{j,k}|^{2N+2}}\leq K,\qquad
(j\geq j_r),
\]
Letting $j\to+\infty$, we obtain
$$\sum_{k} \frac {1-|a_k|^2}{|1-r a_k|^{2N+2}}+\int_0^{2\pi}\frac {d\mu(e^{it})}{|e^{it}-r|^{2N+2}}+\int_0^{2\pi}\frac {|\log|b(e^{it})||}{|e^{it}-r|^{2N+2}}\,dm(e^{it})\leq K$$
for all $r\in (0,1)$. Now we let $r\to 1^-$, we get the desired
condition $\mathrm{(v)}$.

\endpf

\section{Continuity and analytic continuation for functions of the de Branges--Rovnyak spaces}
In this section, we study the continuity and analyticity of
functions in the de Branges--Rovnyak spaces $\HH(b)$ on an open arc
of $\TT$. As we will see the theory bifurcates into two opposite
cases depending whether $b$ is an extreme point of the unit ball of
$H^\infty(\DD)$ or not. Let us recall that if $X$ is a linear space
and $S$ is a convex subset of $X$, then an element $x\in S$ is
called an extreme point of $S$ if it is not a proper convex
combination of any two distinct points  in $S$. Then, it is well
known (see \cite[page 125]{pD70}) that a function $f$ is an extreme
point of the unit ball of $H^\infty(\DD)$ if and only if
$$\int_\TT\log(1-|f(\zeta)|)\,d\zeta=-\infty.$$

The following result is a generalization of results of Helson
\cite{helson} and Sarason \cite{sarason}. The equivalence of
$\mathrm{(i)}$, $\mathrm{(ii)}$ and $\mathrm{(iii)}$ were proved in
\cite[page 42]{sarason} under the
 assumption that $b$ is an extreme point. Our contribution is the
 last two parts. The mere assumption of continuity implies
 analyticity and this observation has interesting applications.

\begin{Thm}\label{thm:prolongement}
Let $b$ be in the unit ball of $H^\infty(\DD)$ and let $I$ be an open arc of $\TT$. Then the following are equivalent:
\begin{enumerate}
\item[$\mathrm{(i)}$] $b$ has an analytic continuation across $I$ and $|b|=1$ on $I$;

\item[$\mathrm{(ii)}$] $I$ is contained in the resolvent set of
$X^*$;

\item[$\mathrm{(iii)}$] any function $f$ in $\HH(b)$ has an analytic continuation across $I$;

\item[$\mathrm{(iv)}$] any function $f$ in $\HH(b)$ has a continuous extension to  $\DD \cup I$;

\item[$\mathrm{(v)}$] $b$ has a  continuous extension to $\DD \cup I$ and $|b|=1$ on $I$.
\end{enumerate}
\end{Thm}

\beginpf
\noindent $\mathrm{(i)}\Longrightarrow \mathrm{(ii)}$: since $|b|=1$
on an open interval, it is clear that $b$ is an extreme point of the
unit ball of $H^\infty(\DD)$. In that case, we know that the
characteristic function of the operator $X^*$ (in the theory of
Sz-Nagy and Foias) is $b$ (see \cite{sarason86}). But then this
theory tells us that $\sigma(X^*)=\sigma(b)$ (see \cite[Theorem
2.3.4., page 102]{nikolski-controle2}).  Therefore, if $b$ has an
analytic continuation across $I$ and $|b|=1$ on $I$, then $I$ is
contained in the complement of $\sigma(b)$ and thus $I$ is contained
in the resolvent set of $X^*$.

$\mathrm{(ii)}\Longrightarrow \mathrm{(iii)}$: for $f\in\HH(b)$, we
have
\[
f(\omega)=\langle f,k_{\omega}^b\rangle_b=\langle f,(Id-\overline\omega X^*)^{-1}k_0^b\rangle_b.
\]
Now if $I$ is contained in the resolvent set of $X^*$,  then the
vector valued function $\omega\longmapsto (Id-\omega
X^*)^{-1}k_0^b$, thought of as an $\HH(b)$-valued function, can be
continued analytically across $I$ and thus the condition
$\mathrm{(iii)}$ follows.

\noindent $\mathrm{(iii)}\Longrightarrow \mathrm{(iv)}$: is clear.

\noindent $\mathrm{(iv)}\Longrightarrow \mathrm{(v)}$: let
$\omega_0\in\DD$ such that $b(\omega_0)\not=0$. Since
$\dfrac{1-\overline{b(\omega_0)}b(z)}{1-\overline{\omega_0}z}$
belongs to $\HH(b)$, it has a continuous extension to $\DD \cup I$.
Therefore $b$  also has a continuous extension to $\DD \cup I$. Now
let $\zeta_0$ be a point of $I$. An application of the principle of
uniform boundedness shows that the functional on $\H(b)$ of
evaluation at $\zeta_0$ is bounded. Let $k_{\zeta_0}^{b}$ denote the
corresponding kernel function. The family $k_\omega^{b}$  tends
weakly to $k_{\zeta_0}^{b}$ as $\omega$ tends to $\zeta_0$ from
$\DD$. Thus, for any $z\in\DD$, we also have
$$\begin{aligned}
k_{\zeta_0}^{b}(z)=&\langle k_{\zeta_0}^{b},k_z^{b}\rangle_{b}
=\lim_{\omega\to \zeta_0}\langle k_\omega^{b},k_z^{b} \rangle_{b}\\
=&\lim_{\omega\to \zeta_0}
\frac{1-\overline{b(\omega)}b(z)}{1-\overline\omega z}
=\frac{1-\overline{b(\zeta_0)}b(z)}{1-\overline{\zeta_0}z}.
\end{aligned}$$
In particular, the function
$\dfrac{1-\overline{b(\zeta_0)}b(z)}{z-\zeta_0}$ is in $H^2(\CC_+)$,
which is possible only if $|b(\zeta_0)|=1$. Hence we get that
$|b|=1$ on $I$.

\noindent $\mathrm{(v)}\Longrightarrow \mathrm{(i)}$: follows from
standard facts based on the Schwarz's reflection principle.

\endpf

As we have seen in the proof of Theorem \ref{thm:prolongement}, one
of the conditions $\mathrm{(i)}-\mathrm{(v)}$ implies that $b$ is an
extreme point of the unit ball of $H^\infty(\DD)$. Thus, the
continuity (or equivalently, the analytic continuation) of $b$
 or of the elements of $\HH(b)$
on the boundary  completely depend on $b$ being an extreme point or
not. If $b$ is not an extreme point of the unit ball of
$H^\infty(\DD)$ and if $I$ is an open arc of $\TT$, then there
exists necessarily a function $f\in\HH(b)$ such that $f$ has not a
continuous extension to $\DD \cup I$. On the opposite case, if $b$
is an extreme point such that $b$ has  continuous extension to $\DD
\cup I$ with $|b|=1$ on $I$, then all the functions $f\in\HH(b)$ are
continuous on $I$ (and even can be continued analytically across
$I$).

Theorem \ref{thm:prolongement}  shows that the de Branges--Rovnyak
spaces $\HH(b)$ have a remarkable property, i.e. continuity on an
open arc of $\TT$ of all functions of $\HH(b)$ is enough to imply
the analyticity of these functions. This property enables us to show
that the result of Dyakonov \cite{Dyak02} concerning the Bernstein's
inequality in the model spaces is sharp in the sense that we
couldn't extend it to all de Branges--Rovnyak spaces. The definition
of de Branges--Rovnyak spaces of the upper half plane is similar to
its counterpart for the unit disc. First, we make precise a little
more the transfer of the unit disc to the upper half plane $\CC_+$.
We consider $\gamma$ the conformal map from $\CC_+$ onto $\DD$
defined by
$$\gamma(z)=\frac {z-i}{z+i},\qquad z\in\CC_+,$$
and we denote by $U$ the (unitary) map from $L^2(\TT)$ onto $L^2(\RR)$ defined by
\begin{eqnarray}\label{eq:disc-demi-plan}
(Uf)(x):=\frac 1{\sqrt\pi}\frac 1{x+i}f\left(\frac{x-i}{x+i}\right),\qquad x\in\RR,f\in L^2(\TT).
\end{eqnarray}
Then it is well known (see \cite[pages 247-248]{nik81}) that $U$
maps $H^2(\DD)$ onto $H^2(\CC_+)$. Moreover, if $\varphi\in
L^\infty(\TT)$, then
\begin{eqnarray}\label{eq:entrelacement-toeplitz}
U T_\varphi=T_{\varphi\circ\gamma}U.
\end{eqnarray}
Now let $b$ be in the unit ball of $H^\infty(\DD)$ and let
$b_1=b\circ\gamma$. Then, using (\ref{eq:entrelacement-toeplitz}),
basic arguments show that $U$ maps unitarily $\HH(b)$ onto
$\HH(b_1)$. Using this unitary transform, we can obviously state the
analogue of Theorem \ref{Thm:derivative-DD} and Theorem
\ref{thm:prolongement} in the upper half plane $\CC_+$.

\begin{Cor}\label{cor:dyakonov}
Let $b_1$ be a point of the unit ball of $H^\infty(\CC_+)$. Then the following are equivalent:
\begin{enumerate}
\item[$\mathrm{(i)}$] the operator $f\longrightarrow f'$ is a bounded operator from $\HH(b_1)$ into $H^2(\CC_+)$;
\item[$\mathrm{(ii)}$] $b_1$ is an inner function and $b'_1\in H^\infty(\CC_+)$.
\end{enumerate}
\end{Cor}

\beginpf
Using \cite[Theorem 1]{Dyak02}, the only thing to prove is that if
$\mathrm{(i)}$ holds, then $b_1$ is inner. But, if for any function
$f$ in $\HH(b_1)$, we have $f'\in H^2(\CC_+)$, then in particular,
 $f$ has a continuous extension to $\CC_+ \cup \RR$. Thus, using
the analogue of Theorem \ref{thm:prolongement} in the upper half
plane, we see that $b_1$ has a continuous extension to  $\CC_+ \cup
\RR$ and $|b_1|=1$ on $\RR$, which means $b_1$ is an inner function.
\endpf

{\bf Acknowledgments:}  The authors deeply thank the anonymous
referee for his/her valuable remarks and suggestions which
 improved the quality of this paper. We also  thank
A. Baranov for several valuable discussions.


\end{document}